\documentclass[11pt]{article}

\usepackage{a4wide,latexsym,amsmath,varioref}

\newcommand{\documentdate}{27 February 2019}

\renewcommand{\Re}{\hbox{I\hskip -2pt R}}
\newcommand{\smallRe}{\hbox{\footnotesize I\hskip -2pt R}}
\newcommand{\req}[1]{(\ref{#1})}
\newcommand{\beqn}[1]{\begin{equation}\label{#1}}
\newcommand{\eeqn}{\end{equation}}
\newcommand{\tim}[1]{\;\; \mbox{#1} \;\;}
\newcommand{\ms}{\;\;\;\;}
\newcommand{\eqdef}{\stackrel{\rm def}{=}}
\newcommand{\numsection}[1]{\section{#1}\setcounter{equation}{0}}

\newcommand{\bpr}{{\bf Proof.} \hspace{1.5mm}}
\newcommand{\epr}{\hfill $\Box$ \vspace*{1em}}
\newcommand{\proof}[1]{
\begin{list}{}{
\setlength{\topsep}{0.0pt}
\setlength{\partopsep}{0.0pt}
\setlength{\leftmargin}{0.025\textwidth}
\setlength{\rightmargin}{0.5\leftmargin}
\setlength{\labelwidth}{0.5\leftmargin}
\setlength{\labelsep}{0.25\leftmargin}}
\item \bpr #1 \epr \noindent
\end{list}}
\newtheorem{theorem}{Theorem}[section]
\newtheorem{lemma}[theorem]{Lemma}

\newcommand{\llem}[2]{\vspace{\baselineskip} 
\noindent\framebox[\textwidth]{\parbox{0.95\textwidth}{
\begin{lemma} \label{#1} \rm #2 \end{lemma} } } \vspace{\baselineskip} }
\newcommand{\lthm}[2]{\vspace{\baselineskip} 
\noindent\framebox[\textwidth]{\parbox{0.95\textwidth}{
\begin{theorem} \label{#1} \rm #2 \end{theorem} } } \vspace{\baselineskip} }
\newcommand{\barc}{\overline{c}} \newcommand{\barf}{\overline{f}}
\newcommand{\barg}{\overline{g}} \newcommand{\barm}{\overline{m}}
\newcommand{\barJ}{\overline{J}}
\newcommand{\calS}{{\cal S}} \newcommand{\calU}{{\cal U}}
\newcommand{\calL}{{\cal L}} \newcommand{\calF}{{\cal F}} 
\newcommand{\sfrac}[2]{{\scriptstyle \frac{#1}{#2}}}
\newcommand{\half}{\sfrac{1}{2}}
\newcommand{\third}{\sfrac{1}{3}}
\newcommand{\threehalves}{\sfrac{3}{2}}
\newcommand{\quarter}{\sfrac{1}{4}}
\newcommand{\ii}[1]{\{1, \ldots, #1 \}}
\newcommand{\comment}[1]{}
\newcounter{algo}[section]
\renewcommand{\thealgo}{\thesection.\arabic{algo}}
\newcommand{\algo}[3]{\refstepcounter{algo}
\begin{center}\begin{figure}[htbp]
\framebox[\textwidth]{
\parbox{0.95\textwidth} {\vspace{\topsep}
{\bf Algorithm \thealgo : #2}\label{#1}\\
\vspace*{-\topsep} \mbox{ }\\
{#3} \vspace{\topsep} }}
\end{figure}\end{center}}
\newcommand{\barell}{\overline{\ell}}
\newcommand{\barDm}{\overline{\Delta m}}
\newcommand{\barDell}{\overline{\Delta \ell}}
\newcommand{\barpsi}{\overline{\psi}}
\newcommand{\barphi}{\overline{\phi}}

\everymath{\displaystyle}

\title{Minimization of nonsmooth nonconvex functions using inexact
  evaluations and its worst-case complexity}

\author{
  S. Gratton
  \thanks{Universit\'e de Toulouse, INP, IRIT, Toulouse, France.
          Email: serge.gratton@enseeiht.fr },
  E. Simon
  \thanks{Universit\'e de Toulouse, INP, IRIT, Toulouse, France.
          Email: ehouarn.simon@enseeiht.fr},
  and Ph. L. Toint
  \thanks{naXys, University of Namur, Namur, Belgium.
          Email: philippe.toint@unamur.be.
          Partially supported by ANR-11-LABX-0040-CIMI within
          the program ANR-11-IDEX-0002-02.}
}
\date{\documentdate}

\begin{document}

\maketitle
\begin{abstract}
An adaptive regularization algorithm using inexact function and derivatives
evaluations is proposed for the solution of composite nonsmooth nonconvex
optimization.  It is shown that this algorithm needs at most
$O(|\log(\epsilon)|\,\epsilon^{-2})$ evaluations of the problem's functions and
their derivatives for finding an $\epsilon$-approximate first-order stationary
point.  This complexity bound therefore generalizes that provided by
[Bellavia, Gurioli, Morini and Toint, 2018] for inexact methods for smooth
nonconvex problems, and is within a factor $|\log(\epsilon)|$ of the optimal
bound known for smooth and nonsmooth nonconvex minimization with exact evaluations.
A practically more restrictive variant of the algorithm with worst-case complexity
$O(|\log(\epsilon)|+\epsilon^{-2})$ is also presented.
\end{abstract}

\vspace*{1cm}

\section{Introduction}

We consider the problem of finding a local minimum of the following composite problem:

\beqn{prob}
\min_{x \in \smallRe^n} \psi(x) = f(x) + h(c(x)),
\eeqn
where $f$ is a (possibly nonconvex) function from $\Re^n$ into $\Re$ whose
gradient is Lipschitz continuous, $c$ is a (possibly nonconvex) function from
$\Re^n$ into $\Re^m$, whose Jacobian is also Lipschitz continuous, and where
$h$ is a convex (possibly nonsmooth) Lipschitz continuous function from
$\Re^m$ into $\Re$.

Such problems occur in a variety of contexts, like LASSO methods in
computational statistics \cite{Tibs96}, Tikhonov regularization of underdetermined
estimation problems \cite{Hans98}, compressed sensing \cite{Dono06},
artificial intelligence \cite{LecuBottBengHaff98}, penalty or projection methods for constrained
optimization \cite{CartGoulToin11a}, reduced-precision deep-learning
\cite{Wangetal18}, image processing \cite{BeckTebo09}, to cite only
a few examples. We also refer the reader to the excellent review in
\cite{LewiWrig16}. In many of these applications, the function $h$ is cheap to
compute\footnote{For example if $h(x)$ is the Euclidean, $\ell_1$ 
or $\ell_\infty$ norm.}, and its Lipschitz constant is known.

Methods to calculate approximate local solutions of the nonconvex problem
\req{prob} have been studied for many years\footnote{We do not consider here
  the abundant literature on the easier convex case, see \cite{DuchRuan18} for
  a recent instance or \cite{BoydVand04} for a general text.}. If $h$ is
differentiable, standard methods include steepest descent,
Levenberg-Morrison-Marquardt quadratic regularization algorithms or
trust-region techniques (see \cite{ConnGoulToin00}). In this case, the
evaluation complexity (that is the number of times the functions $f$ and $c$
need being evaluated for finding an $\epsilon$-approximate first-order point
is proved to be $O(\epsilon^{-2})$ \cite{Nest04,GratSartToin08}.  Moreover,
this order is known to be optimal \cite{CartGoulToin18b}.  If $h$ is
nonsmooth, applicable methods are the proximal gradient method and its
variants \cite{ReddSraPoczSmol19}, as well as the nonsmooth trust-region and
quadratic regularization methods of \cite{CartGoulToin11a} for the nonconvex
ones. It was also shown in this latter paper that the evaluation compexity
remains $O(\epsilon^{-2})$ despite nonsmoothness.  To the best of the authors'
knowledge, available analysis for nonconvex composite problems requires that
$f$ and $c$ (and often their derivatives) are computed exactly.

Inexact function evaluations are however quite commonly met in practice.  For
instance, $f(x)$ or $c(x)$ may be the result of some truncated iterative
process, making the accuracy of the computed values dependent on the
truncation level. Or $f(x)$ or $c(x)$ could be computed as statistical
estimates (e.g. in subsampling methods for additive problems in machine
learning). Or they may result from the need (for embarked processors) or
desire (for high-end supercomputers) to perform their
evaluation in restricted arithmetic precision whenever possible.
Convergence analysis results for methods with inexact function and/or
derivatives values exist
\cite{Cart93,ConnGoulToin00,XuRoosMaho17,CartSche17,BellGuriMori18,
  LiuLiuHsieTao18,ChenJianLinZhan18,BlanCartMeniSche16,BergDiouKungRoye18,
  BellGuriMoriToin18}
and their practical performance considered \cite{Cart93,GratToin18b},
but all these contributions assume smoothness of the objective function.

The contribution of the present paper is threefold.
\begin{itemize}
\item We first propose a new regularization method for the nonsmooth problem
    \req{prob} that uses dynamic accuracy.
\item We then show that the optimal $O(\epsilon^{-2})$ evaluation complexity
    bound is preserved when using this algorithm, up to a (typically modest)
    factor $|\log(\epsilon)|$.
\item We finally present a variant of the algorithm for which a better complexity
  bound of $O(|\log(\epsilon)|+\epsilon^2)$ can be proved at the price of
  loosing some practicality.
\end{itemize}

\noindent
Our presentation is organized as follows. Section~\ref{algo-s} discusses the
nature of the inexact evaluations and presents the new regularization
algorithm, whose essential properties are then developed in
Section~\ref{props-s}. The corresponding evaluation complexity bound is
derived in Section~\ref{compl-s}. A practically more restrictive variant of
the algorithm with better worst-case complexity is presented in
Section~\ref{variant-s}. What can happen if accuracy is limited is
discussed in Section~\ref{discuss-s}. Conclusions are outlined in
Section~\ref{concl-s}.

\numsection{The Adaptive Regularization Algorithm using Dynamic\\ Accuracy}
\label{algo-s}

As indicated above, we assume that $f$, $c$ and their derivatives are
computed inexactly but that $h$ is exact and its cost negligible compared to
that of obtaining approximate value for $f$, $c$ or their derivatives.
Moreover, we assume that, for some \emph{known} constant $L_h \geq 0$, 
\beqn{h-Lip}
\|h(v)-h(w)\| \leq L_h \|v-w\|
\tim{ for all} v,w\in \Re^m.
\eeqn
where $\|\cdot\|$ denotes the standard Euclidean norm.

Our algorithm is iterative and of the adaptive regularization type. It
constructs a sequence of iterates $\{x_k\}$, at which the function $\psi$ and
its derivatives are computed inexactly. For exact values, we use the notations
\[
f_k \eqdef f(x_k), \ms g_k \eqdef g(x_k) = \nabla_x^1f(x_k), \ms
c_k \eqdef c(x_k), \ms J_k \eqdef J(x_k) = \nabla_x^1c(x_k)
\]
and $ \psi_k \eqdef \psi(x_k)=f_k+h(c_k)$. The ``linearization''
\[
\ell_k(s) \eqdef f_k+g_k^Ts + h(c_k+J_ks). 
\]
will play in important role in what follows.  In particular, we use the fact
that 
\[
\min_{\|d\|\leq 1}\ell_k(d) = \ell_k(0)=\psi_k
\]
if $x_k$ is a local minimizer of \req{prob} \cite[Lemma~2.1]{Yuan85a} to say that $x_k$ is an
$\epsilon$-approximate minimizer if
\beqn{first-order-exact}
\phi_k \leq \epsilon,
\eeqn
where
\beqn{phi-def}
\phi_k \eqdef \phi(x_k)
\eqdef \ell_k(0) - \min_{\|d\|\leq 1}\ell_k(d)
= \max_{\|d\|\leq 1}\Delta \ell_k(d),
\eeqn
with
\[
\Delta\ell_k(s)
\eqdef \ell_k(0)-\ell_k(s)
= -g_k^Ts+h(c_k)-h(c_k+J_ks).
\]
If $x_k$ is not such a point, the standard exact regularization algorithm
\cite{CartGoulToin11a} computes a trial step $s_k$ by approximately
minimizing the regularized model
\beqn{m-def}
m_k(s) \eqdef \ell_k(s) + \frac{\sigma_k}{2}\|s\|^2
\eeqn
over all $s\in \Re^n$, where $\sigma_k$ is an adaptive ``regularization parameter''.
This yields the model decrease $\Delta m_k(s_k)$, where
\[
\Delta m_k(s)
\eqdef m_k(0)-m_k(s)
= -g_k^Ts+h(c_k)-h(c_k+J_ks)-\frac{\sigma_k}{2}\|s\|^2.
\]
The value of the objective function is then computed at the trial point
$x_k+s_k$, which is accepted as the new iterate if the achieved reduction
$\psi_k-\psi(x_k+s_k)$ compares well with the predicted decrease
$\Delta\ell_k(s_k)$.  The regularization parameter $\sigma_k$ is then updated
to reflect the quality of this prediction and a new iteration started.

In our context of inexact values for $f$ and $c$, we will keep the same
general algorithm outline, but will also need to take action to handle the
absolute errors in $f$, $g$, $c$ and $J$, denoted by $\varepsilon_f$,
$\varepsilon_g$, $\varepsilon_c$ and $\varepsilon_J$, respectively. In what
follows, we will \emph{denote inexactly computed quantities with an
overbar}. We assume that absolute errors are bounded, and that approximate
values of $f$, $g$, $c$ and $J$ can be computed, given $\varepsilon_f$,
$\varepsilon_g$, $\varepsilon_c$ and $\varepsilon_J$, as
\beqn{barf-def}
\barf_k \eqdef \barf(x_k,\varepsilon_f)
\tim{ with }
|\barf_k - f_k| \leq \varepsilon_f,
\eeqn
\beqn{barg-def}
\barg_k \eqdef \barg(x_k,\varepsilon_g)
\tim{ with }
\|\barg_k - g(x_k)\| \leq \varepsilon_g,
\eeqn
\beqn{barc-def}
\barc_k \eqdef \barc(x_k,\varepsilon_c)
\tim{ with }
\|\barc_k - c(x_k)\| \leq \varepsilon_c,
\eeqn
\beqn{barJ-def}
\barJ_k \eqdef \barJ(x_k,\varepsilon_J)
\tim{ with }
\|\barJ_k - J(x_k)\| \leq \varepsilon_J.
\eeqn
The accuracy level on $\barf$, $\barg$, $\barc$ and $\barJ$ is thus \emph{dynamic},
in the sense that it is specified by the algorithm in order to ensure its
meaningful progress. We will then consider the inexact objective function
$\barpsi(x) = \barf(x)+h(\barc(x))$, together with its ``linearization'' and
model given by 
\[
\barell_k(s) \eqdef \barf_k+ \barg_k^Ts + h(\barc_k+\barJ_ks)
\tim{and}
\barm_k(s) \eqdef \barell_k(s)+ \frac{\sigma_k}{2}\|s\|^2,
\]
defining their corresponding decreases by
\beqn{barDm-def}
\barDm_k(s)
\eqdef-\barg_k^Ts+h(\barc_k)-h(\barc_k+\barJ_ks)-\frac{\sigma_k}{2}\|s\|^2
\eeqn
and
\beqn{barDell-def}
\barDell_k(s_k)
\eqdef -\barg_k^Ts+h(\barc_k)-h(\barc_k+\barJ_ks).
\eeqn
Finally, the criticality measure $\phi_k$ will be approximated by
\[
\barphi_k
\eqdef \barell_k(0) - \min_{\|d\|\leq 1}\barell_k(d)
= \max_{\|d\| \leq 1}\barDell_k(d).
\]

\noindent
Armed with these definitions, we may establish the following crucial error
bounds.

\llem{accuracy-conds-l}{We have that, for any $k$,
\vspace*{-2mm}
\beqn{err-psi}
|\barpsi_k-\psi_k| \leq \varepsilon_f+ L_h\varepsilon_c 
\eeqn
and, for any $v\in\Re^n$,
\beqn{err-Dell}
| \barDm_k(v) - \Delta m_k(v) |
= | \barDell_k(v) - \Delta \ell_k(v) |
\leq (\varepsilon_g + L_h\varepsilon_J)\|v\| + 2 L_h \varepsilon_c.
\eeqn
}
\proof{
Using successively \req{prob}, the triangle inequality, \req{h-Lip},
\req{barf-def} and \req{barc-def}, we obtain that
\[
\begin{array}{lcl}
  |\barpsi_k-\psi_k|
  &  = & | \barf_k + h(\barc_k) - f_k - h(c_k)|\\*[1ex]
  & \leq & |\barf_k -f_k| + |h(\barc_k)-h(c_k)|\\*[1ex]
  & \leq & \varepsilon_f + L_h\|\barc_k-c_k\|\\*[1ex]
  & \leq & \varepsilon_f + L_h\varepsilon_c
\end{array}
\]
and hence \req{err-psi} holds. Similarly, using now \req{barDm-def},
\req{barDell-def}, the triangle and Cauchy-Schwarz inequalities, \req{h-Lip},
\req{barg-def}, \req{barc-def} and \req{barJ-def}, we deduce that
\[
\begin{array}{lcl}
| \barDm_k(v) - \Delta m_k(v) | 
&   =  & | \barDell_k(v) - \Delta \ell_k(v) | \\
& \leq & |(\barg_k-g_k)^Tv| + | h(\barc_k)-h(c_k)|+ |h(\barc_k+\barJ_kv)- h(c_k+J_kv)|\\*[1ex]
& \leq & \|\barg_k-g_k\|\,\|v\| + L_h\|\barc_k-c_k\|+ L_h\|\barc_k+\barJ_kv-c_k-J_kv\|\\*[1ex]
& \leq & \|\barg_k-g_k\|\,\|v\| + L_h\|\barc_k-c_k\|+ L_h(\|\barc_k-c_k\|+\|\barJ_k-J_k\|\,\|v\|)\\*[1ex]
& \leq & \varepsilon_g  \|v\|   + L_h\varepsilon_c + L_h(\varepsilon_c+\varepsilon_J\|v\|)\\*[1ex]
&   =  & (\varepsilon_g + L_h\varepsilon_J)\|v\| + 2L_h\varepsilon_c.
\end{array}
\]
}  

\noindent
Broadly inspired by \cite{BellGuriMoriToin18}, we may now state our
inexact adaptive regularization algorithm formally, in two stages.  We first
describe its global framework \vpageref{ARLDA}, delegating the more
complicated questions of verifying optimality and computing the step to more
detailed sub-algorithms to be presented in the second stage.

In the ARLDA \footnote{For Adaptive Regularization with Lipschitz model and
  Dynamic Accuracy.} algorithm \vpageref{ARLDA}, $\varepsilon_f^{\max}$,
$\varepsilon_g^{\max}$, $\varepsilon_c^{\max}$ and $\varepsilon_J^{\max}$
stand for upper bounds on $\varepsilon_f$, $\varepsilon_g$, $\varepsilon_c$
and $\varepsilon_J$, and $\omega_k$ can be viewed as an iteration dependent
relative accuracy level on $f$, $\psi$ and the model decreases.

\algo{ARLDA}{The ARLDA Algorithm}
{
\vspace*{-0.3 cm}
\begin{description}
\item[Step 0: Initialization.]  An initial point $x_0$ and an initial
  regularization parameter $\sigma_0>0$ are given, as well as an accuracy
  level $\epsilon \in (0,1)$.  The constants $\alpha$, $\kappa_\omega$,
  $\eta_1$, $\eta_2$, $\gamma_1$, $\gamma_2$, $\gamma_3$,
  $\varepsilon_f^{\max}$, $\varepsilon_g^{\max}$, $\varepsilon_c^{\max}$,
  $\varepsilon_J^{\max}$, $\gamma_\varepsilon$ and $\sigma_{\min}$ are also
  given and satisfy $\sigma_{\min} \in (0, \sigma_0]$,
  \[
  0 < \eta_1 \leq \eta_2 < 1, \;\;
  0< \gamma_1 < 1 < \gamma_2 < \gamma_3, \; \;
  \alpha \in (0,1), \;\; \gamma_\varepsilon \in (0,1).
  \]
  Choose
  $\varepsilon_f \leq \varepsilon_f^{\max}$,
  $\varepsilon_g \leq \varepsilon_g^{\max}$,
  $\varepsilon_c \leq \varepsilon_c^{\max}$,
  $\varepsilon_J \leq \varepsilon_J^{\max}$ and
  $\kappa_\omega \in (0,\third \alpha \eta_1]$ such that
  $
  \omega_0 = \varepsilon_f + L_h\varepsilon_c
  \leq \min[\kappa_\omega,\sigma_0^{-1}].
  $
  Set $k=0$.

\item[Step 1: Compute the optimality measure and check for termination. ]\mbox{}

  If unavailable, compute $\barf_k$, $\barg_k$, $\barc_k$ and $\barJ_k$
  satisfying \req{barf-def}--\req{barJ-def}.
  Apply Algorithm~\ref{step1-a} to check for termination with the iterate
  $x_k$ and $\barpsi(x_k)=\barf_k+h(\barc_k)$, or to obtain
  $\barphi_k> \epsilon/(1+\omega_k)$ if termination does not occur.

\item[Step 2: Step calculation. ]

  Apply Algorithm~\ref{step2-a} to approximately minimize $\barm_k(s)$
  and obtain a step $s_k$ and the corresponding linearized decrease $\barDell_k(s_k)$ such that
  \beqn{suff-decr}
  \barDell_k(s_k)
  \geq \frac{1}{4} \min\left\{1,\frac{\barphi_k}{\sigma_k}\right\}\barphi_k.
  \eeqn

\item[Step 3: Acceptance of the trial point. ] 

  Possibly reduce $\varepsilon_f$ to ensure that
  \beqn{vepsf}
  \varepsilon_f\leq \omega_k\barDell_k(s_k).
  \eeqn
  If $\varepsilon_f$ has been reduced, recompute
  $\barf_k(x_k,\varepsilon_f)$ to ensure \req{barf-def}.
  Then compute $\barf_k(x_k+s_k,\varepsilon_f)$ such that 
  \vspace*{-2mm}
  \beqn{Df+-DT-first}
  |\barf_k(x_k+s_k,\varepsilon_f)-f(x_k+s_k)| \leq \varepsilon_f,
  \eeqn
  set
  $\barpsi(x_k+s_k)= \barf_k(x_k+s_k,\varepsilon_f)+h(\barc_k(x_k+s_k,\varepsilon_c))$,
  $\barpsi_k       = \barf_k(x_k,\varepsilon_f)+h(\barc_k(x_k,\varepsilon_c))$
  and define
  \vspace*{-2mm}
  \beqn{rho-def}
  \rho_k = \frac{\barpsi_k - \barpsi(x_k+s_k)}{\barDell_k(s_k)}.
  \eeqn
If $\rho_k \geq \eta_1$, then define
$x_{k+1} = x_k + s_k$; otherwise define $x_{k+1} = x_k$.

\item[Step 4: Regularization parameter update. ]
Set
\beqn{sigupdate-first}
\sigma_{k+1} \in \left\{ \begin{array}{ll}
{}[\max(\sigma_{\min}, \gamma_1\sigma_k), \sigma_k ]  & \tim{if} \rho_k \geq \eta_2, \\
{}[\sigma_k, \gamma_2 \sigma_k ]          &\tim{if} \rho_k \in [\eta_1,\eta_2),\\
{}[\gamma_2 \sigma_k, \gamma_3 \sigma_k ] & \tim{if} \rho_k < \eta_1.
  \end{array} \right.
\eeqn

\item[Step 5: Relative accuracy update. ]
Set
\vspace*{-2mm}
\beqn{new-acc-first}
\omega_{k+1} = \min \left[ \kappa_\omega,\frac{1}{\sigma_{k+1}}\right]
\eeqn
and redefine
  $\varepsilon_f \leq \varepsilon_f^{\max}$,
  $\varepsilon_g \leq \varepsilon_g^{\max}$,
  $\varepsilon_c \leq \varepsilon_c^{\max}$ and
  $\varepsilon_J \leq \varepsilon_J^{\max}$ 
such that $\varepsilon_f + L_h\varepsilon_c \leq \omega_{k+1}$.
Increment $k$ by one and go to Step~1.
\end{description}
}

\noindent
A few comments on this first view of the algorithm are now useful.
\begin{enumerate}
\item The words ``If unavailable'' at the beginning of Step~1 will turn out to
  be fairly important. In a context where the values of $\barf_k$, $\barg_k$,
  $\barc_k$ and $\barJ_k$ may need to be computed several times but with
  different accuracy requirements in the course of the same iteration (as we
  will see below), they indicate that if one of these function has already
  been computed at the current iterate with the desired accuracy, it need not
  (of course) be recomputed. This imposes the minor task of keeping track of
  the smallest value of the relevant $\varepsilon$ for which each of these
  functions has been evaluated at the current iterate.
\item Observe that the relative accuracy threshold $\omega_k$ is recurred from
  iteration to iteration, and the absolute accuracy requirements
  $\varepsilon_f$, $\varepsilon_g$, $\varepsilon_c$ and $\varepsilon_J$ are
  then determined to enforce the relative error (see \req{vepsf} and
  \req{Df+-DT-first} for $\barf$ at $x_k+s_k$). 
\item However, the redefinition of the absolute accuracy requirements in Step~5
  leaves much freedom.  One possible implementation of this redefinition would
  be to set
 \beqn{vareps-upd}
 \begin{array}{ll}
 \varepsilon_f = \min\Big[\varepsilon_f^{\max},\frac{\omega_{k+1}}{\omega_k}\varepsilon_f\Big],&
 \varepsilon_g = \min\Big[\varepsilon_g^{\max},\frac{\omega_{k+1}}{\omega_k}\varepsilon_g\Big],\\*[1.2ex]
 \varepsilon_c = \min\Big[\varepsilon_c^{\max},\frac{\omega_{k+1}}{\omega_k}\varepsilon_c\Big],&
 \varepsilon_J = \min\Big[\varepsilon_J^{\max},\frac{\omega_{k+1}}{\omega_k}\varepsilon_J\Big],
 \end{array}
 \eeqn
 but this is by no means the only possible choice. In particular, any choice
 of $\varepsilon_g \leq \varepsilon_g^{\max}$ and $\varepsilon_J \leq
 \varepsilon_J^{\max}$ is permitted. Observe that since the sequence
 $\{\sigma_k\}$ produced by \req{sigupdate-first}  (or \req{vareps-upd}) need
 not be monotonically increasing, the sequence $\{\omega_k\}$ constructed in
 \req{new-acc-first} need not be decreasing. We present an alternative to this
 choice in Section~\ref{variant-s}.
\item We will verify in Lemma~\ref{step-ok-l} below that the
  sufficient-decrease requirement \req{suff-decr} is fairly loose.  In fact the constant
  $\quarter$ in this condition can be replaced by any constant in $(0,\half)$
  and/or $\barphi_k$ replaced by $\epsilon$ without affecting our theoretical
  results.
\item When exact functions values can be computed (i.e. $\varepsilon_f =
  \varepsilon_g = \varepsilon_c =  \varepsilon_J = 0$), the ARLDA algorithm
  essentially reduces to the regularization algorithm of
  \cite{CartGoulToin11a}. It is also close in spirit to the AR$p$DA algorithm
  for $p=1$ (AR1DA) of \cite{BellGuriMoriToin18} when $h=0$ and the problem
  becomes smooth, but the step computation is simpler in this reference
  because $\barDell_k(s)$ only involves derivatives' values in that case.
\end{enumerate}

\noindent
The purpose of Step~1 of the ARLDA algorithm is to check for termination by
computing a value of $\barphi_k$ which is \emph{relatively} sufficiently
accurate. As can be expected, computing a relatively accurate value when
$\barphi_k$ itself tends to zero may be too demanding, but we nevertheless
design a mechanism that will allow us to prove (in Lemma~\ref{step1-l} below)
that true $\epsilon$-optimality can be reached in this case.  The details of
the resulting Step~1 are given in Algorithm~\ref{step1-a} \vpageref{step1-a}.
Observe that this algorithm introduces a \emph{possible loop on the accuracy
requirement, between~Step~1.3 and Step~1}.

\algo{step1-a}{Check for termination in Algorithm~\ref{ARLDA}}
{
\vspace*{-0.3 cm}
\begin{description}
\item[Step~1.1. ] Solve
    \beqn{min-phi}
    \max_{\|d\|\leq 1} \barDell_k(d)
    \eeqn
    to obtain a global maximizer $d_k$ and the corresponding $\barDell_k(d_k)$.
\item[Step~1.2.] \mbox{}\\*[-3ex]
    \begin{itemize}
    \item If
       \beqn{phi-rel-ok}
       \varepsilon_g+ L_h\varepsilon_J+ 2L_h\varepsilon_c \leq \omega_k \barDell_k(d_k),
       \eeqn
       then
       \begin{itemize}
          \item define $\barphi_k = \barDell_k(d_k)$;
          \item if $\barphi_k \leq \epsilon/(1+\omega_k)$, terminate the ARLDA
                algorithm with {\tt exit} = 1;
          \item else go to Step~2 of the ARLDA algorithm.
       \end{itemize}
    \item If 
       \beqn{small-Dell}
       \barDell_k(d_k) \leq \half\epsilon
       \tim{and}
       \varepsilon_g+ L_h\varepsilon_J+ 2L_h\varepsilon_c \leq \half \epsilon,
       \eeqn
       terminate the ARLDA algorithm with {\tt exit} = 2.
    \end{itemize}
\item[Step~1.3: ]
    Multiply $\varepsilon_g$, $\varepsilon_c$ and $\varepsilon_J$ by
    $\gamma_\varepsilon$ and restart Step~1 of the ARLDA algorithm.
\end{description}
}

\noindent
Once the algorithm has determined in Step~1 that termination cannot occur at
the current iterate, it next computes $s_k$ in Step~2.  In this computation,
the relative accuracy of the ``linearized decrease'' $\barDell_k(s_k)$ must
again be assessed.  This is achieved in Algorithm~\ref{step2-a}
\vpageref{step2-a}.

\algo{step2-a}{Compute the step $s_k$ in Algorithm~\ref{ARLDA}}
{
\vspace*{-0.3 cm}
\begin{description}
  \item[Step~2.1: ] Solve
    \beqn{min-model}
    \min_{s\in \Re^n} \barm_k(s)
    \eeqn
    to obtain a step $s_k$ together with $\barDm_k(s_k)$ and
    $\barDell_k(s_k)$.
  \item[Step~2.2: ]
    If 
    \beqn{rel-s-ok}
    (\varepsilon_g+ L_h\varepsilon_J)\|s_k\| + 2L_h \varepsilon_c
    \leq \omega_k \barDell_k(s_k),
    \eeqn
    go to Step~3 of the ARLDA algorithm.
\item[Step~2.3: ]
    Otherwise multiply $\varepsilon_g$, $\varepsilon_c$ and $\varepsilon_J$ by
    $\gamma_\varepsilon$ and return to Step~1 of the ARLDA algorithm.
\end{description}
}

\noindent
As for Algorithm~\ref{step1-a}, this algorithm introduces a \emph{possible loop on
the accuracy requirement, between~Step 2.3 and Step~1}. We will show (in
Lemma~\ref{single-iteration-complexity-l} below) that these loops are finite,
and thus that the ARLDA algorithm is well-defined. 

\numsection{Properties of the ARLDA algorithm}\label{props-s}

Having defined the algorithm, we turn to establishing some of its properties,
which will be central to the forthcoming complexity analysis. We first verify
that the requirement \req{suff-decr} can always be achieved.

\llem{step-ok-l}{
A step $s_k$ satisfying 
\beqn{m-decr}
\barDell_k(s_k)
\geq \barDm_k(s_k)
\geq \frac{1}{2}\min\left\{1,\frac{\barphi_k}{\sigma_k}\right\}\barphi_k
\eeqn
(and hence also satisfying \req{suff-decr}) can always be computed.
}
\proof{
  The first inequality results from \req{m-def}. The second is given
  in \cite[Lemma~2.5]{CartGoulToin11a}, and hinges on the convexity of $h$.
}

\noindent
We next show an alternative lower bound on the linearized decrease, directly
resulting from the model's definition.

\llem{Dell-Dm-l}{
  For all $k\geq 0$, we have that
  \beqn{mdecr-ns2}
  \barDell_k(s_k) \geq \frac{\sigma_k}{2}\|s_k\|^2.
  \eeqn
  Moreover, as long as the algorithm has not terminated,
  \beqn{Dell-decr}
  \barDell_k(s_k)
  \geq \delta_k(\epsilon)
  \eqdef \frac{1}{16}\min\left\{1,\frac{\epsilon}{\sigma_k}\right\}\epsilon.
  \eeqn
}

\proof{
Lemma~\ref{step-ok-l} implies that, for all $k$,
\[
0 < \barDm_k(s_k) = \barDell_k(s_k) -\frac{\sigma_k}{2}\|s_k\|^2
\]
and \req{mdecr-ns2} follows. We also have that, using \req{suff-decr} and the
fact that $\barphi_k > \epsilon/(1+\omega_k)$ if termination does not occur,
\[
\barDell_k(s_k)
\geq \frac{1}{4}\min\left\{1,\frac{\epsilon}{\sigma_k(1+\omega_k)}\right\}
\frac{\epsilon}{1+\omega_k}
\]
and \req{Dell-decr} then result from \req{new-acc-first} and the
inequality $\omega_k\leq \kappa_\omega \leq \third\alpha\eta_1< 1$.
} 

\noindent
Our next step is to prove that, if termination occurs, the current iterate is
a first-order $\epsilon$-approximate minimizer, as requested.

\llem{step1-l}{(Inspired by \cite[Lemma~3.2]{BellGuriMoriToin18})
  If the ARLDA algorithm terminates, then
  \beqn{first-order}
  \phi_k \leq \epsilon
  \eeqn
  and $x_k$ is a first-order approximate necessary minimizer.
}

\proof{
Suppose first that the ARLDA algorithm terminates at iteration $k$ with
{\tt exit} = 1 in Step~1.2. From the mechanism of this step, we have that
\req{phi-rel-ok} holds and thus, for each $d$ with $\|d\|\leq 1$
\[
(\varepsilon_g+ L_h\varepsilon_J)\|d\|+ 2L_h\varepsilon_c
\leq \varepsilon_g+ L_h\varepsilon_J  + 2L_h\varepsilon_c
\leq \omega_k\barDell_k(d_k)
\]
As a consequence, \req{err-Dell} ensures that, for all $d$ with $\|d\|\leq 1$,
\[
|\barDell_k(d) - \Delta\ell_k(d)| \leq \omega_k\barDell_k(d_k).
\]
Hence,
\beqn{step1-l-1}
\begin{array}{lcl}
\Delta \ell_k(d)
& \leq & \barDell_k(d) + |\barDell_k(d)-\Delta\ell_k(d)|\\*[1ex]
& \leq & \barDell_k(d_k) + |\barDell_k(d)-\Delta\ell_k(d)|\\*[1ex]
& \leq & (1+\omega_k) \barDell_k(d_k).
\end{array}
\eeqn
where we have used that $\barDell_k(d)\leq \barDell_k(d_k)$ by definition of
$d_k$ to derive the second inequality.
As a consequence, for all $d$ with $\|d\|\leq 1$,
\[
\max\Big\{0,\Delta \ell_k(d)\Big\}
\leq (1+\omega_k) \barDell_k(d_k)
= (1+\omega_k)\barphi_k
\leq \epsilon,
\]
where we have used the definition of $\barphi_k$ to
obtain the last inequality. The conclusion \req{first-order} then follows from
\req{phi-def}.

Suppose now that the ARLDA algorithm terminates with {\tt exit} = 2 (in Step~1.2). We then
obtain, using the first two inequalities of \req{step1-l-1}, 
\req{small-Dell} and \req{err-Dell}, that, for every $d$ such that $\|d\|\leq 1$,
\[
\begin{array}{lcl}
\Delta\ell_k(d)
& \leq & \barDell_k(d_k) + |\barDell_k(d)-\Delta\ell_k(d)|\\*[1ex]
& \leq & \half \epsilon + (\varepsilon_g+L_h\varepsilon_J)\|d\| + 2L_h\varepsilon_c\\*[1ex]
& \leq & \half \epsilon + \varepsilon_g+L_h\varepsilon_J + 2L_h\varepsilon_c\\*[1ex]
& \leq & \half \epsilon + \half\epsilon = \epsilon,
\end{array}
\]
which, combined with \req{phi-def}, again implies \req{first-order} for this case.
}  

\noindent
We now establish a useful property of Step~2 (Algorithm~\ref{step2-a}).

\llem{sbound-l}{Suppose that, at Step~2.2,
  \beqn{slarge}
  \|s_k\| \geq \theta_k
  \eqdef \frac{1}{ \omega_k \sigma_{\min}}\left[
   \varepsilon_g^{\max} + L_h\varepsilon_J^{\max} + \sqrt{(\varepsilon_g^{\max} + L_h\varepsilon_J^{\max})^2+4L_h\varepsilon_c^{\max}}
   \right].
  \eeqn
  Then \req{rel-s-ok} is satisfied and the branch to Step~3 of the ARLDA
  algorithm is executed.
}
\proof{
  Step~2 of the ARLDA algorithm terminates as soon as \req{rel-s-ok} holds,
  which, in view of \req{mdecr-ns2} is guaranteed whenever $\|s_k\|$ exceeds
  the largest root of
  \[
  (\varepsilon_g + L_h\varepsilon_J)\|s_k\| + 2L_h\varepsilon_c
  = \half \omega_k \sigma_k \|s_k\|^2.
  \]
  given by
  \[
  \frac{1}{ \omega_k \sigma_k}\left[
   \varepsilon_g + L_h\varepsilon_J + \sqrt{(\varepsilon_g + L_h\varepsilon_J)^2+4L_h\varepsilon_c\omega_k\sigma_k}
   \right],
  \]
  which is itself bounded above by $\theta_k$ as defined in \req{slarge}
  because of the inequality $\sigma_k \geq \sigma_{\min}$, \req{new-acc-first} and the fact that 
  $\varepsilon_f \leq \varepsilon_f^{\max}$, $\varepsilon_g \leq
  \varepsilon_g^{\max}$, $\varepsilon_c \leq \varepsilon_c^{\max}$ and
  $\varepsilon_J \leq \varepsilon_J^{\max}$.
  } 

\noindent
It is also necessary (as announced above) to prove that the accuracy loops
within iteration $k$ are finite, and thus that the ARLDA algorithm is
well-defined. We therefore give explicit bounds on the maximum number of these
accuracy loops and the resulting number of evaluations of the problem's
inexact functions.

\llem{single-iteration-complexity-l}{
  Each iteration $k$ of the ARLDA algorithm
  involves at most two evaluations of $\barf$ and at most
  $1+\nu_k(\epsilon)$ evaluations of $\barg$, $\barc$ and $\barJ$, where
  $\nu_k(\epsilon)$, the number of times that the accuracy thresholds
  $\varepsilon_g$, $\varepsilon_c$ and $\varepsilon_J$ have been reduced by
  Steps~1.3 or 2.3 at iteration $k$, satisfies the bound 
  \beqn{nmax-g-J}
   \nu_k(\epsilon) \eqdef
   \frac{|\log\big((\varepsilon_g^{\max}+L_h\varepsilon_J^{\max})\max\{1,\theta_k\}+2L_h\varepsilon_c^{\max}\big)
     -\log\big(\omega_k\min\{\half\epsilon,\delta_k(\epsilon)\}\big)|}
  {|\log(\gamma_\varepsilon)|},
  \eeqn
  and where $\delta_k(\epsilon)$ and $\theta_k$ are defined in \req{Dell-decr} and
  \req{slarge}, respectively.
}

\proof{ In order to prove this result, we have to bound the number
  of times the accuracy-improving loops (Step~1.3--Step~1) and
  (Step~2.3--Step~1) are being executed.

  Observe first that, at the beginning of every iteration, $\varepsilon_g$,
  $\varepsilon_J$ and $\varepsilon_c$ are bounded above by
  $\varepsilon_g^{\max}$, $\varepsilon_J^{\max}$ and $\varepsilon_c^{\max}$,
  respectively.  Morever, the mechanism of Algorithms~\ref{step1-a} and
  \ref{step2-a} ensures that they can only be reduced within these
  algorithms, and that this reduction is obtained by multiplication with
  the constant $\gamma_\varepsilon<1$. Thus, if $i$ is the number of times
  $\varepsilon_g$, $\varepsilon_J$ and $\varepsilon_c$ have been reduced in
  Steps~1.3 or 2.3, then
  \beqn{shrink}
  \varepsilon_g \leq \gamma_\varepsilon^i \varepsilon_g^{\max},
  \ms
  \varepsilon_J \leq \gamma_\varepsilon^i \varepsilon_J^{\max}
  \tim{and}
  \varepsilon_c \leq \gamma_\varepsilon^i \varepsilon_c^{\max}.
  \eeqn

  Consider the loop (Step~1.3--Step~1) and suppose that
  \beqn{small-exp}
  \varepsilon_g+ L_h\varepsilon_J+ 2L_h\varepsilon_c
  \leq \half\omega_k\epsilon.
  \eeqn
  First consider the case where $\barDell_k(d_k) \geq \half \epsilon$.
  Combining this last inequality with \req{small-exp} gives that
  \req{phi-rel-ok} holds and thus the loop (Step~1.3--Step~1) is terminated by
  either exiting the ARLDA algorithm with {\tt exit} = 1 or going to its
  Step~2. Suppose now that \req{small-exp} holds and that
  $\barDell_k(d_k) < \half \epsilon$.  Then \req{small-Dell} holds and the
  loop is terminated by exiting the ARLDA algorithm with {\tt exit} = 2.
  Thus, using \req{shrink} and \req{small-exp}, the loop is not activated if
  $i$ is large enough to ensure that
  \beqn{noloop1}
  \gamma_\varepsilon^i \,
  (\varepsilon_g^{\max}+L_h\varepsilon_J^{\max}+2L_h\varepsilon_c^{\max})
  \leq \half \omega_k \epsilon.
  \eeqn
  
  The situation is similar for the loop (Step~2.3--Step~1):
  the mechanism of Algorithm~\ref{step2-a} ensure that the loop is not
  activated when \req{rel-s-ok} holds. Suppose first that $\|s_k\|$ remains
  below $\theta_k$ (as defined in \req{slarge}) for all iterations of the loop 
  (Step~2.3--Step~1). Then, in view of \req{shrink} and \req{Dell-decr},
  \req{rel-s-ok} must hold at the latest when
  \beqn{noloop2}
  \gamma_\varepsilon^i \,
  \Big((\varepsilon_g^{\max}+L_h\varepsilon_J^{\max})\theta_k+2L_h\varepsilon_c^{\max}\Big)
  \leq \omega_k\delta_k(\epsilon)
  \eeqn
  where $\delta_k(\epsilon)$ is defined in \req{Dell-decr}. If $\|s_k\|$
  happens to exceed $\theta_k$ before \req{noloop2} is satisfied, then
  \req{rel-s-ok} is also satisfied earlier because of Lemma~\ref{sbound-l} and
  the loop terminated. We therefore deduce from \req{noloop1} and
  \req{noloop2} that the loops (Step~1.3--Step~1) and (Step~2.3--Step~1) can
  be activated at most $\nu_k(\epsilon)$ times during the complete $k$-th
  iteration of the ARLDA algorithm, where $\nu_k(\epsilon)$ is given by
  \req{nmax-g-J}. Thus $\barg(x_k,\varepsilon_g)$, $\barc(x_k,\varepsilon_c)$
  and $\barJ(x_k,\varepsilon_j)$ are computed (in the beginning of Step~1) at
  most $1+\nu_k(\epsilon)$ times.  The observation that $\barf_k$ is computed
  at most two times per ARLDA iteration (in Step~3) concludes the proof.  }

\noindent
We next bound the error on the successive values of the objective function.

\llem{ps-psi+-err-l}{We have that, for all $k\geq 0$,
\beqn{err-psi-2}
|\barpsi_k - \psi_k| \leq \threehalves \omega_k \barDell_k(s_k)
\tim{ and }
|\barpsi_k^+ - \psi_k^+| \leq \threehalves \omega_k \barDell_k(s_k).
\eeqn
}  

\proof{
When $\rho_k$ is computed in Step~3, it must be because Step~2 has been
completed, and hence \req{rel-s-ok} must hold, which in turn implies that
$L_h\varepsilon_c\leq \half \omega_k \barDell_k(s_k)$.  Thus the desired
inequalities follow from \req{barf-def}, \req{Df+-DT-first}, \req{vepsf} and
\req{err-psi}.
}

\noindent
We finally recall a standard result on successful versus unsuccessful iterations.

\llem{SvsU-l}{
\cite[Theorem~2.4]{BirgGardMartSantToin17}
Let
\vspace*{-2mm}
\beqn{calSk-def}
\calS_k = \{ j \in \ii{k} \mid \rho_j \geq \eta_1\}
\tim{ and }
\calU_k = \ii{k}\setminus\calS_k
\eeqn
be the sets of \emph{successful} and \emph{unsuccessful} iterations, respectively.
The mechanism of Algorithm~\ref{ARLDA} guarantees that, if
\beqn{sigmax}
\sigma_{k} \leq \sigma_{\max},
\vspace*{-2mm}
\eeqn
for some $\sigma_{\max} > 0$, then
\vspace*{-2mm}
\beqn{unsucc-neg}
k +1 \leq |\calS_k| \left(1+\frac{|\log\gamma_1|}{\log\gamma_2}\right)+
\frac{1}{\log\gamma_2}\log\left(\frac{\sigma_{\max}}{\sigma_0}\right).
\eeqn
}

\noindent
This shows that it is sufficient, for establishing the overall evaluation
complexity of the ARLDA algorithm, to bound the maximum number of
evaluations at successful iterations.

\numsection{Worst-case evaluation complexity}\label{compl-s}

We are now in position to start our evaluation complexity proper.  In order to
make it formally coherent, we start by explicitly stating our assumptions on
the problem.

\begin{description}
\item[AS.1.] $f$ and $c$ are continuously differentiable in $\Re^n$.
\item[AS.2.] There exist non-negative constants 
  $L_g$ and $L_J$ such that, for all $k\geq 0$,
  and for all $x,y$ in
  $\calL_0 = \{ v \in \Re^n \mid \psi(v) \leq \psi(x_0) \}$,
  \beqn{gJ-Lip}
  \|g(x)-g(y)\| \leq 2L_g \|x-y\|
  \tim{and}
  \|J(x)-J(y)\| \leq 2L_J \|x-y\|.
  \eeqn
\item[AS.3] There exists a constant $L_h\geq 0$ such that \req{h-Lip} holds.
\item[AS.4] There exists a constant $\psi_{\rm low}$ such that $\psi(x)\geq
  \psi_{\rm low}$ for all $x \in \Re^n$.
\end{description}

\noindent
A first (and standard) consequence of AS.1-AS.2 is the following result on
error bounds for $f$ and $c$ at a trial point $x+s$.

\llem{Lip-ders}{Suppose that AS.1 and AS.2 hold.  Then, for all $x,s\in \Re^n$,
\[
|f(x+s) - (f(x) + g(x)^Ts )| \leq L_g\|s\|^2
\tim{and}
\|c(x+s) - (c(x) + J(x)s )\| \leq L_J\|s\|^2
\]
}

\proof{See \cite[Lemma~2.1]{CartGoulToin18b}.} 

\noindent
We may then use the bounds to establish the following important bound on the
regularization parameter.

\llem{sig-bounded-l}{Suppose that AS.1-AS.3 hold. Then there exists a constant
  $\sigma_{\max} \geq \max\{1,\sigma_0\}$ such that, for all $k\geq 0$, 
\beqn{sig-bound}
\sigma_k \leq \sigma_{\max}
\eqdef \max\left\{\sigma_0,\gamma_3\frac{4 + 2(L_g + L_hL_J)}{1-\eta_2},\frac{1}{\kappa_\omega}\right\}
\tim{and}
\omega_k \geq \frac{1}{\sigma_{\max}}.
\eeqn
}

\proof{
We have that
\[
\begin{array}{lcl}
|\rho_k-1|
&  =   & \frac{|\barpsi_k-\barpsi_k^+-\barDell_k+\Delta\ell_k - \Delta\ell_k|}{\barDell_k}\\
& \leq & \frac{1}{\barDell_k}\Big[ | \barpsi_k-\psi_k| +
  |\barpsi_k^+-\psi_k^+| + |\barDell_k - \Delta\ell_k| +
  |\psi_k^+ -(\psi_k-\Delta\ell_k)|\Big]\\
& \leq & \frac{1}{\barDell_k}\Big[ 4\omega_k\barDell_k
  + |\psi_k^+ -(\psi_k-\Delta\ell_k)|\Big],
\end{array}
\]
where we also used \req{rho-def}, the triangle inequality to derive the first
inequality, while the second results from \req{err-psi-2} and the fact that,
if the algorithm has not terminated at iteration $k$, then \req{rel-s-ok} must
hold, in turn implying \req{err-Dell} because of Lemma~\ref{accuracy-conds-l}.
Now, because of the triangle inequality, \req{gJ-Lip}, Lemma~\ref{Lip-ders}
and \req{h-Lip}, we see that
\[
\begin{array}{lcl}
|\psi_k^+ -( \psi_k-\Delta\ell_k)|
&   =  & | f_k^+ + h(c_k^+)- f_k -h(c_k)  -g_k^Ts_k+ h(c_k) - h(c_k+J_ks_k)| \\*[1ex]
& \leq & | f_k^+-(f_k+g_k^Ts_k)| + |h(c_k^+)-h(c_k+J_ks_k)|\\*[1ex]
& \leq & | f_k^+-(f_k+g_k^Ts_k)| + L_h\|c_k^+-c_k+J_ks_k\|\\*[1ex]
& \leq & L_g\|s_k\|^2 + L_hL_J\|s_k\|^2.
\end{array}
\]
Thus, combining the two last displays,
\beqn{sig-l-1}
|\rho_k-1|
\leq 4 \omega_k+ (L_g + L_hL_J)\frac{\|s_k\|^2}{\barDell_k}.
\eeqn
Taking now \req{new-acc-first} and the inequality of \req{mdecr-ns2}
into account, we deduce that
\[
|\rho_k-1|
\leq \frac{1}{\sigma_k}\Big[ 4 + 2(L_g + L_hL_J)\Big]
\leq 1-\eta_2
\tim{ whenever }
\sigma_k \geq \frac{4 + 2(L_g + L_hL_J)}{1-\eta_2},
\]
in which case $\rho_k\geq \eta_2 \geq \eta_1$, iteration $k$ is successful
(i.e. $k \in \calS_k$) and $\sigma_{k+1}\leq \sigma_k$. The mechanism of the
algorithm then ensures that \req{sig-bound} holds for all $k$.
The lower bound on $\omega_k$ follows from \req{new-acc-first} and the fact
that \req{sig-bound} ensures that $(1/\sigma_{\max}) \leq \kappa_{\omega}$.
} 

\noindent
The bound \req{sig-bound} is important, in particular because it allows, in
conjunction with AS.2, to simplify the bound on the complexity of a single
iteration of the ARLDA algorithm, making this bound only dependent on $\epsilon$
(i.e. dropping the dependence on $k$).

\llem{better-single-l}{Suppose that AS.1-AS.3 hold.  Then we have that, before
  termination, each iteration of the ARLDA algorithms evaluates $\barf$ at
  most two times and $\barc$, $\barg$ and $\barJ$ at most $1+\nu(\epsilon)$
  times, where
  \beqn{nu-def}
  \nu(\epsilon)
  \eqdef   \frac{|2\log\big(\epsilon\big)|
    + |\log\big((\varepsilon_g^{\max}+L_h\varepsilon_J^{\max})\theta +2L_h\varepsilon_c^{\max}\big)
    + 2\log(4\sigma_{\max})
    |}
  {|\log(\gamma_\varepsilon)|}
  \eeqn
  with
  \beqn{theta-def}
  \theta \eqdef \max\left\{1,
  \frac{\sigma_{\max}}{\sigma_{\min}}\left[
   \varepsilon_g^{\max} + L_h\varepsilon_J^{\max} + \sqrt{(\varepsilon_g^{\max} + L_h\varepsilon_J^{\max})^2+4L_h\varepsilon_c^{\max}}
   \right]
  \right\}.
  \eeqn
}

\proof{
  We observe that, because of \req{Dell-decr}, \req{sig-bound} and the
  inequalities $\epsilon \leq 1 \leq \sigma_{\max}$ and the second part of \req{sig-bound},
  \[
  \omega_k \delta_k(\epsilon)
  \geq
  \frac{\omega_k}{16}\min\left\{1,\frac{\epsilon}{\sigma_{\max}}\right\}\epsilon
  \geq \frac{\epsilon^2}{16\sigma_{\max}^2}
  \tim{ and }
  \half\omega_k\epsilon \geq \frac{\epsilon^2}{16\sigma_{\max}^2}.
  \]
  Moreover, the second part of \req{sig-bound} and \req{slarge} imply that $\theta_k \leq \theta$,
  a value independent of $k$ and $\epsilon$. Using these bounds in \req{nmax-g-J}, we see
  that
  \[
  \nu_k(\epsilon)
  \leq \frac{\log\left(\frac{\epsilon^2}{16\sigma_{\max}^2}\right)
    -\log\big((\varepsilon_g^{\max}+L_h\varepsilon_J^{\max})\theta +2L_h\varepsilon_c^{\max}\big)}
  {\log(\gamma_\varepsilon)}
  \]
  which, with Lemma~\ref{single-iteration-complexity-l},the second part of
  \req{sig-bound} and the observation that the above value only depends on
  $\epsilon$, concludes the proof. 
}  

\noindent
Following a well-worn path in complexity analysis, we may now use a
telescopic sum argument involving successive objective function's decreases at
successful iterations and Lemmas~\ref{step-ok-l},  \ref{SvsU-l}
and \ref{better-single-l} to deduce our final result.

\lthm{final}{Suppose that AS.1-AS.4 hold.  Then the ARLDA algorithm terminates with
  $\phi_k \leq \epsilon$ in at most
  \[
  \tau(\epsilon)\;\;\mbox{iterations},\;\;
  2\tau(\epsilon) \;\;\mbox{evaluations of $\barf$, and}\;\;
  \lfloor (1 + \nu(\epsilon)) \tau(\epsilon) \rfloor \;\;\mbox{evaluations of $\barg$, $\barc$ and 
   $\barJ$},
  \]
  where
  \beqn{tau-def}
  \tau(\epsilon) \eqdef \left\lfloor
  \frac{8\sigma_{\max}\big(\psi(x_0)-\psi_{\rm low}\big)}{\eta_1(1-\alpha)}\,\epsilon^2 +1\right\rfloor
  \left(1+\frac{|\log\gamma_1|}{\log\gamma_2}\right)+
     \frac{1}{\log\gamma_2}\log\left(\frac{\sigma_{\max}}{\sigma_0}
     \right),
  \eeqn
  $\nu(\epsilon)$ is defined in \req{nu-def} and $\sigma_{\max}$ is
  defined in \req{sig-bound}.
}

\proof{
If iteration $k$ is successful (i.e. $k \in \calS_k$) and the ARLDA algorithm
has not terminated yet, one has that
\[
\begin{array}{lcl}
\psi(x_k)-\psi(x_{k+1})
& \geq & [\barpsi_k(x_k) - \barpsi_k(x_{k+1})]
         - 3 \omega_k\barDell_k(s_k) \\*[2ex]
& \geq & \eta_1 \barDell_k(s_k)
         -\alpha \eta_1\barDell_k(s_k)\\*[2ex]
& \geq &\frac{\eta_1
   (1-\alpha)}{2}\min\left\{1,\frac{\barphi_k}{\sigma_k}\right\} \barphi_k\\*[2ex]
& \geq &\frac{\eta_1
   (1-\alpha)}{2}\min\left\{1,\frac{\epsilon}{\sigma_{\max}(1+\omega_k)}\right\} \frac{\epsilon}{1+\omega_k}\\*[2ex]
& = &\frac{\eta_1(1-\alpha)\epsilon^2}{2\sigma_{\max}(1+\omega_k)^2},
\end{array}
\]
where we used \req{err-psi-2}, \req{rho-def}, \req{m-decr} and \req{sig-bound},
the fact that $\barphi_k>\epsilon/(1+\omega_k)$ before termination, that
$\sigma_{\max}\geq 1$ and the inequality $\epsilon\leq 1$. Thus $\psi(x_k)$ is
monotonically decreasing, and one then deduces that
\[
\psi(x_0)-\psi(x_{k+1})
\geq \frac{\eta_1(1-\alpha)\epsilon^2}{2\sigma_{\max}(1+\omega_k)^2} \,|\calS_k|.
\]
Using that $\psi$ is bounded below by $\psi_{\rm low}$ and the inequalities
$\omega_k\leq\kappa_\omega<1$, we conclude that 
\[
| \calS_k |
\leq \frac{2\sigma_{\max}(1+\omega_k)^2}{\eta_1(1-\alpha)}(\psi(x_0) - \psi_{\rm low}) \epsilon^{-2}
< \frac{8\sigma_{\max}}{\eta_1(1-\alpha)}(\psi(x_0) - \psi_{\rm low}) \epsilon^{-2}
\]
until termination. Lemmas~\ref{SvsU-l} and \ref{sig-bounded-l} are then
invoked to compute the upper bound on the total number of iterations
$\tau(\epsilon)$, and Lemma~\ref{better-single-l} is invoked to bound the
number of evaluations.
}

\noindent
If, as is usual in evaluation complexity analysis, one focuses on the maximum
number of evaluations expressed as the order in $\epsilon$, the bound of
Theorem~\ref{final} may be simplified to
\beqn{in-order}
O\Big( |\log(\epsilon)|\, \epsilon^{-2} \Big) \tim{evaluations,}
\eeqn
which is identical in order to the bound obtained for the inexact first-order
regularization method AR1DA in \cite{BellGuriMoriToin18}.

\numsection{An algorithmic variant with monotonic accuracy thresholds}\label{variant-s}

As in \cite{BellGuriMoriToin18}, we now consider a variant of the ARLDA
algorithm for which a better worst-case complexity bound can be proved, at
the price of a signifiucantly more rigid dynamic accuracy strategy.

Suppose that the relatively loose conditions for updating
$\varepsilon_f$, $\varepsilon_g$, $\varepsilon_c$ and $\varepsilon_J$ and the
end of Step~5 of the ARLDA algorithm are replaced by
\beqn{variant}
\mbox{If necessary, \fbox{decrease} $\varepsilon_f$, $\varepsilon_g$,
  $\varepsilon_c$ and $\varepsilon_J$ to ensure that
  $\varepsilon_f+L_h\varepsilon_c \leq \omega_{k+1}$.}
\eeqn

In this case, $\varepsilon_f$, $\varepsilon_g$, $\varepsilon_c$ and
$\varepsilon_J$ all decrease monotonically. As a consequence, the number of
times they are reduced by multiplication with $\gamma_\epsilon$ is still
bounded by $\nu(\epsilon)$ as given in \req{nu-def}, but this bound now holds
for reductions at Steps~1.3 or 2.3 \emph{across all iterations} (instead of at
iteration $k$ only). We may therefore revise Theorem~\req{final} as follows.

\lthm{final-variant}{Suppose that AS.1-AS.4 hold.  Then the variant
  of the ARLDA algorithm using the update \req{variant} terminates with
  $\phi_k \leq \epsilon$ in at most
  \[
  \tau(\epsilon)\;\;\mbox{iterations},\;\;
  2\tau(\epsilon) \;\;\mbox{evaluations of $\barf$, and}\;\;
  \lfloor \nu(\epsilon) + \tau(\epsilon) \rfloor \;\;\mbox{evaluations of $\barg$, $\barc$ and 
   $\barJ$},
  \]
  where $\tau(\epsilon)$ is defined in \req{tau-def} and
  $\nu(\epsilon)$ is defined in \req{nu-def}.
}

\proof{
The proof is identical to that of Theorem~\ref{final} except for the very
last argument, where one now needs to take the revised interpretation of $\nu(\epsilon)$
into account to derive the maximum number of approximate evaluations of $g$, $c$
and $J$.
} 

\noindent
Observe that expressing this new bound in order of $\epsilon$ now gives
\[
O\Big(|\log(\epsilon)| + \epsilon^{-2}\Big) \tim{evaluations,}
\]
which typically improves upon \req{in-order} and extends the bound known in
the smooth case for the $p=1$ variant of the AR$p$DA algorithm
with monotonic accuracy \cite{BellGuriMoriToin18}.
But, as indicated above this improved bound comes at
the price of the more restrictive updating rule \req{variant}. In particular
this rule means that a (potentially large) number of iterations will require
an accuracy on $g$, $c$ and $J$ which is tighter than what is actually needed
for the algorithm's progress.

\numsection{Discussion}\label{discuss-s}

The theory presented above supposes a somewhat ideal world, where arbitrarily
high accuracy may be requested for the evaluation of the problem's function
values and their derivatives. In practice however, such requests are likely to
be too demanding, for instance due to limitations of computer arithmetic. It
may thus happen that evaluating $\barf$, $\barc$, $\barg$ or $\barJ$ becomes
impossible, especially if $\psi$ is locally very nonlinear causing $\sigma_k$
to increase and $\omega_k$ to decrease.

A first comment is that algorithmic precautions may be taken, in the framework
of the present theory, to make this event less likely. The most obvious one
is to use the ARLDA algorithm itself (instead of its variant of
Section~\ref{variant-s}). Secondly, it is important to choose the final
accuracy $\epsilon$ large enough to ensure that satisfying
\beqn{safe0}
\varepsilon_g+L_h\varepsilon_J +2L_h\varepsilon_c
\approx \half \epsilon
\eeqn
(the second inequality in \req{small-Dell}) is at all possible. Moreover, as
one expects $\barDell(s_k)$ to be of the order of $\epsilon$ and $\|s_k\|$ to
be of the order of $\sqrt{\epsilon}$ when converging, \req{rel-s-ok} and
\req{new-acc-first} suggest that the condition 
\beqn{safe1}
(\varepsilon_g+L_h\varepsilon_J)\sqrt{\epsilon} +2L_h\varepsilon_c
\approx \frac{\epsilon}{\sigma_{\max}}
\eeqn
should be achievable, where $\sigma_{\max}$ is given by \req{sig-bound}.
Assuming the term in $\sigma_0$ does not dominate is this latter expression,
the condition \req{safe1} becomes
\beqn{safe1b}
(\varepsilon_g+L_h\varepsilon_J)\sqrt{\epsilon} +2L_h\varepsilon_c
\approx \frac{\epsilon(1-\eta_2)}{\gamma_3(3+2(L_g + L_hL_J))}.
\eeqn
Similarly, \req{vepsf} and \req{new-acc-first} indicate that
\beqn{safe2}
\varepsilon_f \approx \frac{\epsilon(1-\eta_2)}{\gamma_3(3+2(L_g + L_hL_J))}
\eeqn
should also be achievable.  This discussion furthermore indicates that
limiting the growth of $\sigma_k$ as much as possible by choosing moderate
values of $\gamma_2$ and $\gamma_3$ in \req{sigupdate-first} might be a good
idea. A third possibility is to ``balance'' the accuracy requests between
$\varepsilon_g$, $\varepsilon_J$ and $\varepsilon_c$ in order to satisfy
\req{phi-rel-ok} and \req{rel-s-ok}, depending on the value of $L_h$. For
instance, if $L_h$ is large, one might consider choosing $\varepsilon_g$
smaller to allow for a larger $\varepsilon_c$.  In view of \req{rel-s-ok},
this is even more important if $\|s_k\|$ is small (as can be expected when
converging). Finally, since \req{vepsf} and \req{rel-s-ok} involve
$\barDell_k(s_k)$ in their right-hand side, computing the step $s_k$ more
accurately than requested by \req{suff-decr} may also be helpful.

As indicated, these stategies may still be insufficient because the high
nonlinearity inherent to the problem causes $\sigma_k$ to grow or because the
conditions \req{safe0}--\req{safe2} are too restrictive to hold in practice.
If failure to compute one of the problem's function occurs with values of
$\sigma_k$ barely ensuring successful iterations, we contend that this is
signal that the algorithm should be stopped as it has exhausted its ``descent
potential'' on the exact objective function.  Three cases must be considered.
The first is when the value of $\barDell_k(d_k)$ cannot be proved to be
significant enough for its value to be interpreted as the optimality measure
$\barphi_k$ (this likely to happen for quite small values of
$\barDell_k(d_k)$). This implies that the link between $\phi_k$ and
$\barphi_k$ is lost, but the proof of Lemma~\ref{step1-l} nevertheless
indicates that ``noisy optimality'' is achieved in the sense that, for all $d$
with $\|d\|\leq 1$,
\[
\Delta \ell_k(d)
\leq \max\{\half \epsilon,\barDell_k(d_k)\} + \varepsilon_g+L_h\varepsilon_J
     + 2L_h\varepsilon_c.
\]
The second case is when \req{rel-s-ok} cannot be satisfied, meaning that
$\barDell_k(s_k)$ cannot be made accurate enough (due to failing evaluations
of $\barc$, $\barg$ or $\barJ$) to make the latter significant compared with
the inaccuracy noise. Because of the form of \req{rel-s-ok}, it is possible
that backtracking along the step $s_k$ could improve the situation, as
convexity of $\ell_k$ leaves the possibility that $\|\beta s_k\|$ decreases
faster that $\barDell_k(\beta s_k)$ for $\beta$ tending to zero in $(0,1]$,
thereby allowing \req{rel-s-ok} to hold for some $\beta$.  If this is the
case, minimization can be pursued, possibly at the price of loosing the
complexity guarantee of Theorem~\req{final} if $\barDell_k(\beta s_k)$ is
too small compared to $\epsilon^2$.  If \req{rel-s-ok} cannot be enforced,
this means that progress based on the model cannot be guaranteed, and the
algorithm should then be stopped.  A similar situation occurs in the third
case, where the computation of $\barf(x_k,\varepsilon_f)$ or
$\barf(x_k+s_k,\varepsilon_f)$ fails.  This then means that the decrease in
the objective function value is obscured by inaccuracies and cannot be
meaningfully compared to the predicted decrease.  A purely deterministic
algorithm, like ARLDA, must therefore abandon. But, as we have noted, it is
not because the correct working of the method is no longer \emph{guaranteed}
that significant objective function decrease may not happen by
chance. Attempting some re-evaluations and/or recomputations of $\barDell_k$
may, with some luck, allow progress.  It is therefore not unreasonable to
consider such an effort-limited ``trial-and-error'' heuristic, close to
random-direction search, if the algorithm stalls due to impractical accuracy
requests. Obviously, this is beyond the theory we have presented.

We conclude this section by an important observation. Since the mechanism of
requiring adaptive absolute errors on the inexactly computed quantities is
identical to that used in \cite{BellGuriMoriToin18}, the probabilistic
complexity analysis derived in this reference remains valid for our case.
Moreover, if either $f$ or $c$ is computed by subsampling sums of many
independent terms (as is frequent in machine learning applications), the
sample size estimators presented in \cite[Theorem~6.2]{BellGuriMoriToin18} may
also be used in our framework.

\numsection{Conclusion and perspectives}\label{concl-s}

For solving the possibly nonsmooth and nonconvex composite problem \req{prob},
we have proposed an adaptive regularization algorithm using inexact
evaluations of the problem's functions and their first derivative, whose
evaluation complexity is $O\big( |\log(\epsilon)|\, \epsilon^{-2}\big)$.  This
complexity bound is within a factor $|\log(\epsilon)|$ of the known optimal
bound for first-order methods using exact derivatives for smooth
\cite{CartGoulToin18b} or nonsmooth composite \cite{CartGoulToin11a}
problems. It also generalizes the bound derived in \cite{BellGuriMoriToin18}
to the composite nonsmooth case. We have also shown that a practically more
restrictive variant of the algorithm has
$O\big(|\log(\epsilon)|+\epsilon^{-2}\big)$ complexity.

Our method and analysis can easily be extended to cover two other cases of
potential interest. The first is when $g$ and $J$ are merely
$\beta$-H\"{o}lder continuous rather that Lipschitz-continuous, and the second
is to set-constrained problems $\min\psi(x)$ for $x\in\calF$, where the
constraints are \emph{inexpensive} in the sense that
their/evaluation/enforcement has a negligible cost compared to that of
evaluating $f$, $g$, $c$ or $J$. We have refrained from including the
generality needed to cover these two extensions here for clarity of
exposition, and we refer the reader to
\cite{CartGoulToin18b,BellGuriMoriToin18} for details. We also note that, as
in \cite{CartGoulToin18b} (for instance), the Lipschitz conditions of AS.2
need only to apply on each segment of the ``path of iterates''
$\cup_{k\geq0}[x_k,x_k+1]$ for our results to hold.

The authors are aware that there is considerable room for an updating
strategy for $\varepsilon_f$, $\varepsilon_g$, $\varepsilon_c$ and
$\varepsilon_J$ which is more practical than uniform multiplication
by $\gamma_\varepsilon$ or simple updates of the form \req{vareps-upd}.
One expects their worst-case complexity to lie between $O\big( |\log(\epsilon)|\,
\epsilon^{-2}\big)$ and $O\big(|\log(\epsilon)|+\epsilon^{-2}\big)$ depending
on how much non-monotonicity is allowed. They should be considered in a
(desirable) numerical study of the new methods. 

{\footnotesize

\section*{Acknowledgment}

The third author is grateful to the ENSEEIHT (INP, Toulouse) for providing a
friendly research environment for the duration of this research project. 
  

\begin{thebibliography}{10}

\bibitem{BeckTebo09}
A.~Beck and M.~Teboulle.
\newblock A fast iterative shrinkage-thresholding algorithm for linear inverse
  problems.
\newblock {\em SIAM J. Imaging Sci.}, 2:183--202, 2009.

\bibitem{BellGuriMori18}
S.~Bellavia, G.~Gurioli, and B.~Morini.
\newblock Theoretical study of an adaptive cubic regularization method with
  dynamic inexact {H}essian information.
\newblock arXiv:1808.06239, 2018.

\bibitem{BellGuriMoriToin18}
S.~Bellavia, G.~Gurioli, B.~Morini, and {Ph.}~L. Toint.
\newblock Deterministic and stochastic inexact regularization algorithms for
  nonconvex optimization with optimal complexity.
\newblock arXiv:1811.03831, 2018.

\bibitem{BergDiouKungRoye18}
E.~Bergou, Y.~Diouane, V.~Kungurtsev, and C.~W. Royer.
\newblock A subsampling line-search method with second-order results.
\newblock arXiv:1810.07211, 2018.

\bibitem{BirgGardMartSantToin17}
E.~G. Birgin, J.~L. Gardenghi, J.~M. Mart\'{i}nez, S.~A. Santos, and Ph.~L.
  Toint.
\newblock Worst-case evaluation complexity for unconstrained nonlinear
  optimization using high-order regularized models.
\newblock {\em Mathematical Programming, Series~A}, 163(1):359--368, 2017.

\bibitem{BlanCartMeniSche16}
J.~Blanchet, C.~Cartis, M.~Menickelly, and K.~Scheinberg.
\newblock Convergence rate analysis of a stochastic trust region method via
  supermartingales.
\newblock arXiv:1609.07428v3, 2018.

\bibitem{BoydVand04}
S.~Boyd and L.~Vandenberghe.
\newblock {\em Convex Optimization}.
\newblock Cambridge University Press, Cambridge, England, 2004.

\bibitem{Cart93}
R.~G. Carter.
\newblock Numerical experience with a class of algorithms for nonlinear
  optimization using inexact function and gradient information.
\newblock {\em SIAM Journal on Scientific and Statistical Computing},
  14(2):368--388, 1993.

\bibitem{CartGoulToin11a}
C.~Cartis, N.~I.~M. Gould, and Ph.~L. Toint.
\newblock On the evaluation complexity of composite function minimization with
  applications to nonconvex nonlinear programming.
\newblock {\em SIAM Journal on Optimization}, 21(4):1721--1739, 2011.

\bibitem{CartGoulToin12a}
C.~Cartis, N.~I.~M. Gould, and Ph.~L. Toint.
\newblock On the oracle complexity of first-order and derivative-free
  algorithms for smooth nonconvex minimization.
\newblock {\em SIAM Journal on Optimization}, 22(1):66--86, 2012.

\bibitem{CartGoulToin18b}
C.~Cartis, N.~I.~M. Gould, and Ph.~L. Toint.
\newblock Sharp worst-case evaluation complexity bounds for arbitrary-order
  nonconvex optimization with inexpensive constraints.
\newblock arXiv:1811.01220, 2018.

\bibitem{CartSche17}
C.~Cartis and K.~Scheinberg.
\newblock Global convergence rate analysis of unconstrained optimization
  methods based on probabilistic models.
\newblock {\em Mathematical Programming, Series~A}, 159(2):337--375, 2018.

\bibitem{ChenJianLinZhan18}
X.~Chen, B.~Jiang, T.~Lin, and S.~Zhang.
\newblock On adaptive cubic regularization {N}ewton's methods for convex
  optimization via random sampling.
\newblock arXiv:1802.05426, 2018.

\bibitem{ConnGoulToin00}
A.~R. Conn, N.~I.~M. Gould, and Ph.~L. Toint.
\newblock {\em Trust-Region Methods}.
\newblock MPS-SIAM Series on Optimization. SIAM, Philadelphia, USA, 2000.

\bibitem{Dono06}
D.~L. Donoho.
\newblock Compressed sensing.
\newblock {\em IEEE Trans. Inform. Theory}, 52(4):1289--1306, 2006.

\bibitem{DuchRuan18}
J.~Duchi and F.~Ruan.
\newblock Stochastic methods for composite and weakly convex optimization
  problems.
\newblock {\em SIAM Journal on Optimization}, 28(4):3229--3259, 2018.

\bibitem{GratSartToin08}
S.~Gratton, A.~Sartenaer, and Ph.~L. Toint.
\newblock Recursive trust-region methods for multiscale nonlinear optimization.
\newblock {\em SIAM Journal on Optimization}, 19(1):414--444, 2008.

\bibitem{GratToin18b}
S.~Gratton and Ph.~L. Toint.
\newblock A note on solving nonlinear optimization problems in variable
  precision.
\newblock arXiv:1812.03467, 2018.

\bibitem{Hans98}
P.~C. Hansen.
\newblock {\em Rank-Deficient and Discrete Ill-Posed Problems: Numerical
  Aspects of Linear Inversion}.
\newblock SIAM, Philadelphia, USA, 1998.

\bibitem{LecuBottBengHaff98}
Y.~LeCun, L.~Bottou, Y.~Bengio, and P.~Haffner.
\newblock Gradient-based learning applied to document recognition.
\newblock {\em Proceedings of the IEEE}, 86(11):2278--2324, 1998.

\bibitem{LewiWrig16}
A.~S. Lewis and S.~J. Wright.
\newblock A proximal method for composite minimization.
\newblock {\em Mathematical Programming, Series~A}, 158:501--546, 2016.

\bibitem{LiuLiuHsieTao18}
L.~Liu, X.~Liu, C.-J. Hsieh, and D.~Tao.
\newblock Stochastic second-order methods for non-convex optimization with
  inexact {H}essian and gradient.
\newblock arXiv:1809.09853, 2018.

\bibitem{Nest04}
{Yu}. Nesterov.
\newblock {\em Introductory Lectures on Convex Optimization}.
\newblock Applied Optimization. Kluwer Academic Publishers, Dordrecht, The
  Netherlands, 2004.

\bibitem{ReddSraPoczSmol19}
S.~Reddi, S.~Sra, B.~P\'{o}czos, and A.~Smola.
\newblock Proximal stochastic methods for nonsmooth nonconvex finite-sum
  optimization.
\newblock In D.~D. Lee, M.~Sugiyama, U.~V. Luxburg, I.~Guyon, and R.~Garnett,
  editors, {\em Advances in Neural Information Processing Systems 29}, pages
  1145--1153. Curran Associates, Inc., 2016.

\bibitem{Tibs96}
R.~Tibshirani.
\newblock Regression shrinkage and selection via the {LASSO}.
\newblock {\em Journal of the Royal Statistical Society B}, 58(1):267--288,
  1996.

\bibitem{Wangetal18}
N.~Wang, J.~Choi, D.~Brand, C.-Y. Chen, and K. Gopalakrishnan.
\newblock Training deep neural networks with 8-bit floating point numbers.
\newblock In {\em 32nd Conference on Neural Information Processing Systems},
  2018.

\bibitem{XuRoosMaho17}
P.~Xu, F.~Roosta-Khorasani, and M.~W. Mahoney.
\newblock {N}ewton-type methods for non-convex optimization under inexact
  {H}essian information.
\newblock arXiv:1708.07164v3, 2017.

\bibitem{Yuan85a}
Y.~Yuan.
\newblock Conditions for convergence of trust region algorithms for nonsmooth
  optimization.
\newblock {\em Mathematical Programming}, 31(2):220--228, 1985.

\end{thebibliography}

}

\end{document}